# Do CRM Systems Cause One-to-One Marketing Effectiveness?

**Sunil Mithas, Daniel Almirall and M. S. Krishnan**

*Abstract.* This article provides an assessment of the causal effect of customer relationship management (CRM) applications on one-to-one marketing effectiveness. We use a potential outcomes based propensity score approach to assess this causal effect. We find that firms using CRM systems have greater levels of one-to-one marketing effectiveness. We discuss the strengths and challenges of using the propensity score approach to design and execute CRM related observational studies. We also discuss the applicability of the framework in this paper to study typical causal questions in business and electronic commerce research at the firm, individual and economy levels, and to clarify the assumptions that researchers must make to infer causality from observational data.

*Key words and phrases:* Causal analysis, potential outcomes, propensity score, matching estimator, customer relationship management systems, electronic commerce.

## 1. INTRODUCTION

Electronic commerce, information systems and marketing researchers widely agree on the importance of causal analysis to gain a better understanding of the causal effects of managerial interventions and marketing programs (Boulding, Staelin, Ehret and Johnston, 2005; Gregor, 2006; Lucas, 1975). Does implementing an information technology (IT) intervention, such as a customer relationship management (CRM) system, improve firm performance (Brynjolfsson and Hitt, 1998; Mithas, Krishnan and Fornell, 2005; Srinivasan and Moorman, 2005)? Do managerial interventions, such as customer satisfaction improvement programs, add to shareholder wealth (Fornell, Mithas, Morgeson and Krishnan, 2006; Peppers, Rogers and Dorf, 1999; Rust and Kannan, 2003)? Do certain auction parameters maximize seller or consumer surplus (Bapna, Goes and Gupta, 2001; Bapna, Jank and Shmueli, 2005; Koppius and Van Heck, 2002; Mithas and Jones, 2006)? Does acquiring an MBA degree cause an increase in the salary of an IT professional (Connolly, 2003; Mithas and Krishnan, 2004b)? Has "offshoring" caused a decline in the jobs and wages in the U.S. economy (Mithas and Whitaker, 2006; Venkatraman, 2004)? The common theme across these questions—at the firm, individual and economy levels—is that they all are causal questions posed by managers, individuals and policy makers about the causal effects of some intervention.

Historically, with their early emphasis on individual level phenomena, such as the impact of mode of information presentation on learning and performance (Lucas and Nielsen, 1980), electronic commerce researchers answered causal questions using an experimental approach that permitted randomization—the gold standard for assessing causality. As researchers increasingly focus on firm level phenomena, particularly the business value of IT interventions, they must rely on observational data, because randomized field trials across actual firms are virtually impossible. In addition, the growing digitization of business processes is providing unprecedented opportunities to collect richer observational data (more for individual level behavior than for firm level research), making it easier to pose new questions that researchers did not or could not ask before (Jank and Shmueli, 2006). However, the use of observational (nonexperimental) data in these settings raises concerns regarding the ability to interpret causal results from empirical analyses (Rosenbaum, 1999). What exactly is the problem with using observational data? What benefits does randomization provide in the experimental setting and what assumptions are needed to make causal claims using nonexperimental data? Can statistical methods help answer questions of causality? The potential outcomes framework for causal inference, described in more detail by Rubin (1974) and Holland (1986), sheds light on these questions and provides researchers with tools to answer some of the questions posed in this Introduction.

*Sunil Mithas is Assistant Professor, Department of Decision and Information Technologies, Robert H. Smith School of Business, University of Maryland, College Park, Maryland 20742, USA e-mail: smithas@umd.edu. Daniel Almirall is a graduate student in Statistics, Institute for Social Research, University of Michigan, Ann Arbor, Michigan 48104, USA e-mail: dalmiral@umich.edu. M. S. Krishnan is Professor and Hallman Fellow, Department of Business Information Technology, Stephen M. Ross School of Business, University of Michigan, Ann Arbor, Michigan 48109, USA e-mail: mskrish@umich.edu.*

In this paper, we point to the usefulness of a potential outcomes approach, called propensity score stratification (Rosenbaum and Rubin, 1983b), to investigate causal relationships in the substantive domain of electronic commerce research. For further discussion of the potential outcomes approach, refer to Angrist and Krueger (1999), Heckman (2005), Imbens (2004), Rosenbaum (2002), Rubin (2005) and Winship and Morgan (1999). We address a problem in the electronic commerce domain and estimate the causal effect of CRM systems on one-to-one marketing effectiveness. The approach relies on the assumption of strong ignorability, which implies that assignment to a treatment group is independent of potential outcomes conditional on observed pretreatment covariates (Rosenbaum and Rubin, 1983b).

## 2. CRM SYSTEMS AND ONE-TO-ONE MARKETING EFFECTIVENESS

Firms invest over $50 billion each year on IT applications (such as CRM systems) to streamline customer-interfacing business processes. A primary objective of these systems is to improve one-to-one marketing effectiveness, that is, the ability of a firm to target an individual customer based on previous history and purchasing behavior. However, media reports question whether these CRM implementations have paid off (Harvard Management Update, 2000), and from a business value perspective (Banker, Kauffman and Mahmood, 1993; Barua and Mukhopadhyay, 2000; Brynjolfsson and Hitt, 1998; Kauffman and Weill, 1989; Sambamurthy, Bharadwaj and Grover, 2003), there is a need to estimate whether CRM systems have indeed caused an improvement in one-to-one marketing effectiveness.

We obtained data for this study from *InformationWeek* magazine. The data include information about firms' IT systems and related business benefits for the year 1999, and were collected by *InformationWeek* between late 1999 and early 2000 as part of a more comprehensive survey to benchmark firms' IT infrastructure and managerial practices in their respective industries. *InformationWeek* has been surveying top IT managers (including Vice Presidents, Chief Information Officers and Directors) of large firms in the United States since 1986 to identify the firms that are the best users of IT. *InformationWeek* surveys are typically sent to very large firms such as Fortune 1000 firms. Because *InformationWeek* has a significant presence in the IT business community and offers visibility and incentives to participating firms, we believe that the high response rate for these surveys compares favorably to other firm level academic studies. *InformationWeek* is considered to be a reliable source of information, and previous academic studies have also used data from *InformationWeek* surveys (Mithas, Krishnan and Fornell, 2005; Rai, Patnayakuni and Patnayakuni, 1997). Our sample consists of 487 firms.

We now specify the effect we wish to estimate by using the language of potential outcomes (Rubin, 2005). Let $t$ denote a firm's CRM status: $t = 1$ if a firm adopts CRM and $t = 0$ if a firm does not adopt CRM. Let the potential outcome $Y(t)$ denote a firm's assessment of its one-to-one marketing effectiveness from its customer related IT applications. For a fixed $t$, $Y(t) = 1$ if the firm is effective in one-to-one marketing and $Y(t) = 0$ if a firm is not effective. Observe that each firm has two potential response values: $Y(t = 1)$ is the firm's effectiveness had the firm adopted CRM and $Y(t = 0)$ is the firm's effectiveness had the firm not adopted CRM. Causal effects are defined as contrasts in these two quantities.

Using this notation, we denote the average causal effect $\Delta$ of CRM on $Y$ as the difference in the proportion of firms with effective one-to-one marketing programs had *all firms* adopted CRM versus the proportion with effective programs had *none of the firms* adopted CRM. The average causal effect $\Delta$ describes the change in the proportion of companies with effective marketing programs that is caused by CRM. Note that lowercase $t$ is not a random variable, but is instead conceived as an index of the response $Y(t)$ which is assumed to exist for every firm in our study. In other words, neither $t$ nor $Y(t)$ represents observed data. The observed data analogs of $t$ and $Y(t)$ are represented by the random variables uppercase $Z$ and $Y$, respectively. That is, we denote a firm's observed CRM status using uppercase $Z$. In this study, 282 of the 487 firms adopted a CRM system. The observed outcome variable $Y$ is a firm's assessment of the improvement in one-to-one marketing effectiveness from its customer related IT applications. This observed response is binary (1 indicates that a firm reported an increase in one-to-one marketing effectiveness; 0 indicates that a firm did not report an increase in one-to-one marketing effectiveness). Note that the observed one-to-one marketing effectiveness is only one of the two potential outcomes, because we cannot observe the same firm both with CRM and without CRM.

Finally, the three-vector $X$ is a set of observed pre-CRM characteristics, including (1) MFG—a firm's industry sector (manufacturing or services firm), (2) ITPC—a firm's amount of IT investment as a percentage of revenue and (3) CUSTAPPS—the presence of other customer related IT systems (this 13-item scale indicates deployment of IT systems to support the business processes involved in customer acquisition and disposal of products and services offered by firms, including product marketing information, multilingual communication, personalized marketing offerings, dealer locator, product configuration, price negotiation, personalization, transaction system, on-line distribution and fulfillment system, customer service and customer satisfaction tracking) at the firm prior to CRM adoption. Table 1 shows summary statistics for the treatment (CRM) and control (non-CRM) groups.

A naïve estimate of $\Delta$ is the difference in proportion, which we denote by $D$, of observed CRM and observed non-CRM firms that report increases in one-to-one marketing effectiveness, that is, $D = E(Y|Z=1) - E(Y|Z=0)$, which in our sample is equal to 0.30 [SE = 0.04, 95% CI = (0.22, 0.39)]. However, it is problematic to attribute the causal effect $\Delta$ to the observed difference in the proportion of CRM versus non-CRM firms that report increases in one-to-one marketing effectiveness. This is because we do not observe the performance of firms with CRM *had they not adopted CRM* and vice versa, and because firms that adopt CRM systems may be different in ways related to performance outcomes (i.e., confounding or selection bias) (Rosenbaum and Rubin, 1983b). Thus, we may have a selection problem; that is, in symbols, $D = \Delta + \text{Bias}$, where Bias is the extent to which the CRM and non-CRM groups differ according to pre-treatment variables (observed or unobserved) that also predict $Y$.

In an experimental setting where firms are randomly assigned to either a treatment or control group, these selection problems would not occur. The act of randomizing CRM in an experimental setting, for example, would render treatment (CRM) and control (non-CRM) groups equal and balanced on average (both in observed and unobserved characteristics). That is, Bias = 0 on average. The treatment effect can then be assessed by comparing the mean outcomes of treatment and control groups at a given time after the treatment is assigned. In contrast, in a nonexperimental setting such as actual firms, CRM status may be affected by one or more of a firm's observed and/or unobserved characteristics. For example, the use of CRM systems may be due to self-selection by managers or may be mandated by industry consortiums or business partners (Mithas, Krishnan and Fornell, 2005). Comparing the mean outcomes in observational studies may therefore overestimate or underestimate the true causal effect of CRM implementation, and even a larger sample size would not remedy these selection problems.

TABLE 1
*Characteristics of treatment and control groups before matching*

|  | Non-CRM[a] (control group) $N = 199$ | CRM[a] (treatment group) $N = 282$ | Mean difference (p value) |
|---|---|---|---|
| **Outcome variable** | | | |
| Improvement in one-to-one | 0.37 | 0.67 | 0.30 |
| marketing effectiveness | (0.48) | (0.47) | (<0.01) |
| **Observed covariates** | | | |
| Customer-facing IT systems (CUSTAPPS) | 6.77 | 8.24 | 1.47 |
|  | (2.51) | (2.79) | (<0.01) |
| IT investments (ITPC) | 2.91 | 3.59 | 0.68 |
|  | (2.43) | (3.53) | (<0.01) |
| Manufacturing (MFG) | 0.45 | 0.34 | $-0.11$ |
|  | (0.50) | (0.48) | (<0.01) |

[a]The mean and standard deviation (in parentheses) are given.

## 3. A PROPENSITY SCORE APPROACH TO ESTIMATE CAUSAL EFFECT

We begin by checking the differences between CRM and non-CRM firms on observed pre-CRM covariates. Figure 1 shows the distribution of observed covariates across CRM and non-CRM firms. As Figure 1 and Table 1 show, compared to non-CRM firms, CRM firms have more customer-facing IT systems, invest more in IT as a percentage of revenues and are underrepresented in the manufacturing sector. We conducted a more formal analysis of the selection into CRM status using a probit model with CRM status as the dependent variable and observed covariates as explanatory variables. The chi-square test in the probit model reveals that the selection model is significant compared to a model with no explanatory variables. Thus, CRM firms differ significantly from non-CRM firms with respect to observable covariates in the probit model. Characteristics such as prior investments in customer-facing IT systems and a firm's industry sector significantly affect the probability of adopting CRM. Because it is likely that observed pretreatment covariates are also associated with observed one-to-one marketing effectiveness, this analysis suggests that Bias may be nonzero.

The following estimator of $\Delta$ relies on the assumption of strong ignorability (Rosenbaum and Rubin, 1983b). This assumption says that selection bias is due only to correlation between observed firm characteristics and a firm's treatment status, where treatment refers to the implementation or absence of CRM systems. Formally, the assumption states that the distribution of the observed CRM status does not depend on the potential outcomes $Y(0)$ and $Y(1)$ given the observed covariates. We have already shown that there exists some evidence of selection based on observables. By making the strong ignorability assumption, we are saying that there exists no unmeasured or unknown pre-treatment variable (say $U$) that is correlated directly with both CRM status and the outcome variable one-to-one marketing effectiveness. [A firm's market orientation (Kohli and Jaworski, 1990) is an example of one such unmeasured covariate $U$ that may violate this assumption, because market orientation affects CRM adoption and may also have a direct effect on customer related outcomes such as one-to-one marketing effectiveness. Other unmeasured pre-treatment covariates may exist, including variables not yet known to electronic commerce researchers. A procedure with R code that describes the sensitivity of our results according to violations of this assumption is available from the authors on request.]

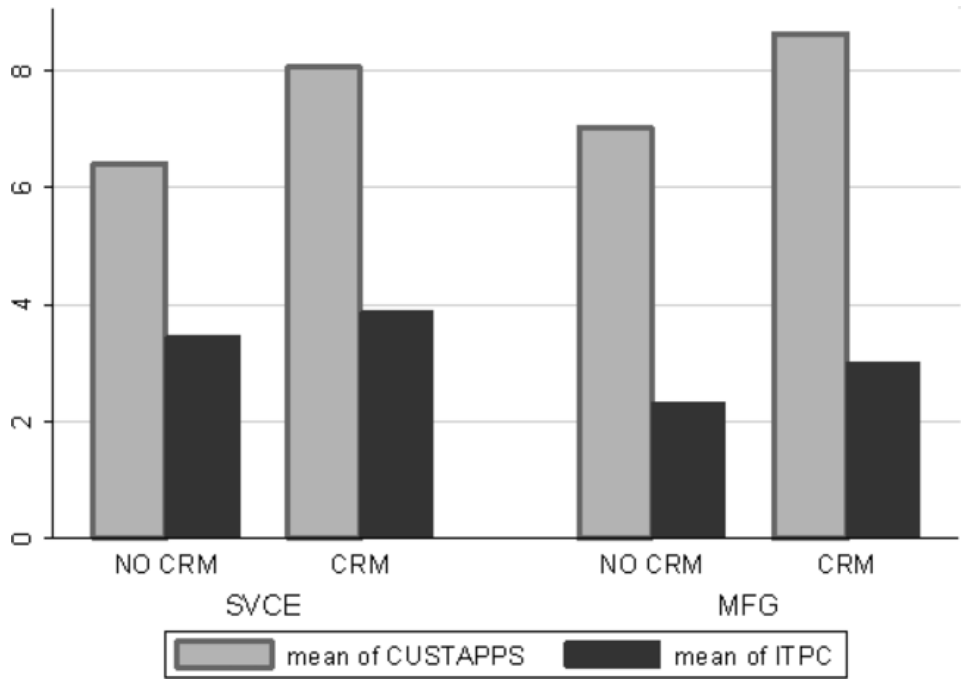

Fig. 1. *Covariate balance before stratification.*

Under the assumption of strong ignorability, matching on the basis of observed covariates can be used to overcome the selection bias problem. With matching, the basic idea is to find a non-CRM firm (i.e., a "clone" in the sense of Rubin and Waterman, 2006) for every CRM firm such that the firms do not differ in any way other than their CRM adoption status. Although not a serious issue in our CRM application that has only three covariates, multivariate matching is usually problematic because sample sizes are often not large enough to achieve matching on all observed covariates. This problem (known as the curse of dimensionality) becomes particularly severe if the covariates are of a continuous nature. Extending the previous work of Rubin (1977) that shows the efficacy of balancing on a single covariate, Rosenbaum and Rubin (1983b) proposed an ingenious solution to this problem. Their idea is to match on the propensity score instead. The propensity score is defined as the probability of being treated given observed covariates—which, in our setting, translates to the probability that a firm has adopted CRM given the covariates. We denote the propensity score by $e(x) = \Pr(Z = 1 | X = x)$.

We apply a matching technique based on the estimated value of propensity scores $e(x)$, a function of the observed pre-CRM characteristics. This approach, also known as propensity score stratification, relies on forming subclasses (or strata) $S(x)$ based on the propensity scores (Rosenbaum and Rubin, 1984). It can be seen as a way to further reduce the dimensionality of observed covariates [i.e., $X$ is reduced to $e(x)$, which is reduced to the subclasses $S(x)$]. Specifically, the objective of subclassification is to create subclasses based on the propensity score so that CRM and non-CRM firms have similar values of the propensity score (thereby achieving balance on the multivariate $X$). This enables a fair comparison (on $Y$) of CRM and non-CRM firms within each subclass. Then, once an appropriate level of balance within each subclass is achieved [e.g., between $e(x) = 0.2$ and $e(x) = 0.4$], an estimate of $\Delta$ can be obtained by taking a weighted average of the within strata effects. Rosenbaum and Rubin (1983b) show that it is sufficient to condition on the univariate propensity score $e(x)$ [and thus $S(x)$] to remove bias due to multivariate $X$. To construct the subclasses, we use the propensity scores based on the probit model mentioned above (Becker and Ichino, 2002). (We also tried interaction and quadratic terms in our propensity score models. Because they did not affect the balance of covariates significantly across the CRM and non-CRM groups, we did not use them in our final analysis.)

Because the matching estimators do not identify the treatment effect outside the region of common support, our first step is to calculate the range of support for both the CRM and non-CRM groups. In our study, the support for the CRM group is (0.22–0.89) and the support for the control group is (0.18–0.91). Following Rubin (2001), we dropped from our analysis six non-CRM firms that fell outside the common support region (propensity scores for these firms are either less than 0.22 or more than 0.89). As Rubin and Waterman (2006) elaborate, it is appropriate to drop these firms from our analysis because for these firms, we cannot find an appropriate "clone." Our second step is to create the subclasses along the space of common support. Following Dehejia and Wahba (2002), we initially classified all observations in five equal-sized subclasses based on propensity scores. Because there are no firms in the 0–0.2 propensity score range, this meant starting with four subclasses (with inferiors of propensity scores at 0.2, 0.4, 0.6 and 0.8). Then we checked for any differences in propensity scores across CRM and non-CRM firms in each stratum. If we found any significant differences, we subdivided the stratum until we obtained a similar distribution of propensity scores and covariates in each stratum. This resulted in five strata that achieved propensity score and covariate balance across CRM and non-CRM firms.

TABLE 2
*Covariate balance after subclassification on propensity score*[a]

| Stratum | Customer-facing IT systems | | IT investments | | Manufacturing | |
|---|---|---|---|---|---|---|
| | CRM | Non-CRM | CRM | Non-CRM | CRM | Non-CRM |
| 1 | 3.65 | 3.64 | 2.02 | 2.14 | 0.73 | 0.71 |
| | (0.30) | (0.23) | (0.24) | (0.23) | (0.09) | (0.08) |
| 2 | 6.17 | 6.12 | 2.70 | 2.61 | 0.41 | 0.45 |
| | (0.16) | (0.15) | (0.24) | (0.18) | (0.05) | (0.05) |
| 3 | 8.79 | 8.81 | 3.50 | 3.34 | 0.36 | 0.42 |
| | (0.19) | (0.23) | (0.39) | (0.51) | (0.06) | (0.08) |
| 4 | 10.86 | 10.36 | 4.84 | 4.55 | 0.18 | 0.14 |
| | (0.14) | (0.35) | (0.46) | (0.68) | (0.04) | (0.07) |

[a]The mean and standard deviation (in parentheses) are given for each category. Note that none of the covariates has a statistically significant difference across the CRM and non-CRM firms within a stratum.

Because of our limited sample size, the fifth stratum had only four non-CRM firms. Therefore, we combined the fourth and fifth strata to have a reasonable number of CRM and non-CRM units in each stratum. Thus, our final analysis contained four subclasses. (We also tried three and six subclasses. However, our attempts to use three and six subclasses did not succeed because this prevented a balancing of covariates in each stratum, one of the critical conditions for subsequent analysis. Using more than six subclasses is not practical in our case because of the limited sample size that is typical in firm level research, unlike the common use of nine or ten subclasses in individual level studies that have several thousand observations.) Table 2 shows the covariate balance and summary statistics across CRM and non-CRM firms within each stratum after subclassification based on propensity scores. We note that the covariate balance after propensity score stratification is better (see Figure 2 and Table 2) com-pared to that before stratification (see Figure 1 and Table 1).

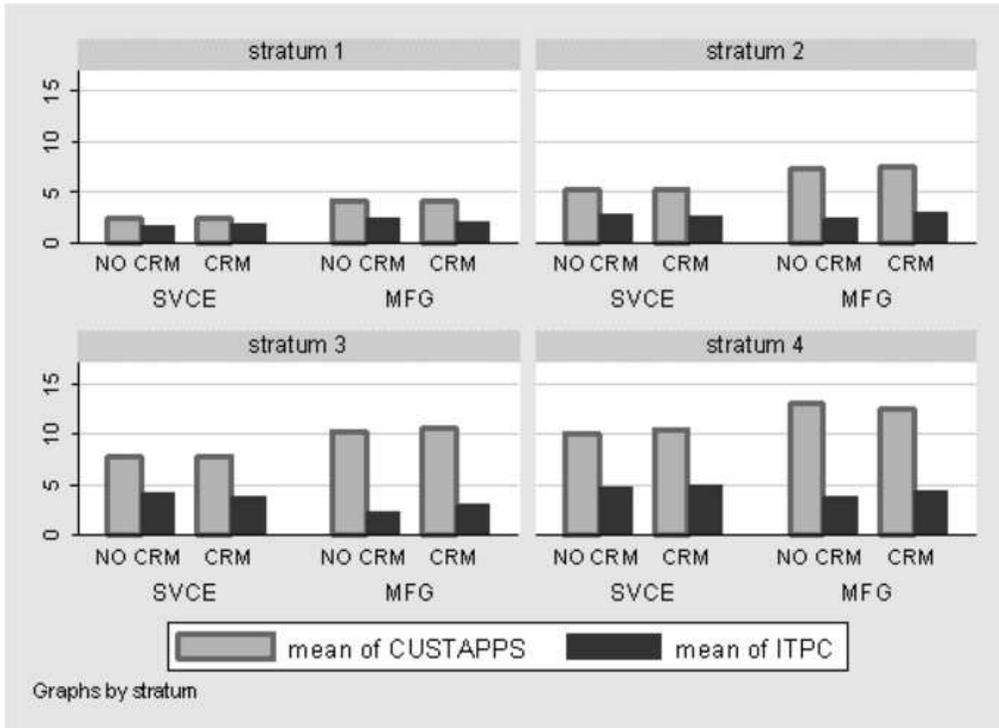

FIG. 2. *Covariate balance after stratification.*

TABLE 3
*Rubin's diagnostics for assessing covariate balance before and after stratification*

| | $B$ Standardized difference in means | $R1$ Ratio of variances of propensity score | $R2$, variance ratio orthogonal to the propensity score | | |
|---|---|---|---|---|---|
| | | | Manufacturing | Customer-facing IT systems | IT investments |
| Before stratification | 0.65 | 1.34 | 0.91 | 0.93 | 1.97 |
| After four subclass stratification | 0.22 | 1.16 | 1.05 | 0.96 | 1.25 |

We used diagnostics suggested by Rubin (2001) to further verify the success of the propensity score stratification in achieving covariate balance. We computed the following: (1) $B$, the standardized difference in propensity score means in the CRM and non-CRM groups (recommended value less than $1/2$), (2) $R1$, the ratio of the propensity score variances in the two groups (recommended value between 0.8 and 1.25) and (3) $R2$, the ratio of variances of the residuals of each covariate (MFG, CUSTAPPS and ITPC) after adjusting for the propensity score (recommended value between 0.80 and 1.25). Table 3 presents the results of these calculations and shows the efficacy of four subgroups based on propensity score stratification in achieving the covariate balance, because the values of $B$, $R1$ and $R2$ fall within the recommended range after stratification. In particular, stratification helped us bring $B$, $R1$ and $R2$ closer to recommended values.

Table 4 shows the average effect of CRM on marketing effectiveness within each stratum. The average effect (over the four subclasses) of CRM on one-to-one marketing is 23 percentage points. This effect takes into account selection bias due to correlation between three observed variables (used in the selection equation) and the treatment variable. This analysis suggests that the naïve estimator $D$ overestimated $\Delta$ in this application. One reason for this may be that the presence of a greater number of customer-facing IT systems is positively associated with both use of CRM systems and one-to-one marketing effectiveness, causing Bias to be positive.

We now compare the results obtained using the propensity score approach with those obtained using a probit model. If there is a substantial difference in the composition of treatment and control groups, classical regression adjustment methods may not be reliable. Regression adjustment methods make strong functional form assumptions (e.g., constant treatment effect) based on extrapolation, and thus may not provide an estimate of the causal effect one hopes to estimate. Despite this consideration, researchers who wish to employ classical regression methods (such as ordinary least squares, log-linear models or logistic regression) to adjust for pretreatment covariates will find the propensity score $e(x)$ to be a useful diagnostic tool that provides better visibility of the extent to which treatment and control groups differ. For example, Rubin's diagnostics in row 1 of Table 3 (before stratification) suggest that in our CRM application, regression adjustment may reliably adjust for selection bias because the values of $B$, $R1$ and $R2$ are not significantly different from the recommended values (except for the value of $R2$ for IT investments, all other values are only slightly outside the recommended range). Therefore, we expect to see a similar estimated treatment effect if we use a regression technique. To confirm this, we ran a probit model of one-to-one marketing effectiveness on all covariates used in the propensity score model and the CRM indicator variable. We find that CRM applications are positively associated with an improvement in one-to-one marketing effectiveness ($\beta = 0.59$, $p < 0.001$). Because the coefficients of probit models are not easily interpretable, we also discuss the effect of a unit change in the CRM variable on the probability of gain in one-to-one marketing effectiveness. We find that CRM systems are associated with a 23.37% increase in the probability of an increase in one-to-one marketing effectiveness, holding all other variables constant at their mean values. The results of this analysis suggest the usefulness of the propensity score approach as a first step to test for covariate balance before making subsequent regression adjustments.

TABLE 4
*Causal effect of CRM*

| Stratum | CRM treated units | Non-CRM control units | Effect of CRM on one-to-one marketing effectiveness (standard errors) |
|---|---|---|---|
| 1 | 23 | 31 | 0.35 (0.13) |
| 2 | 90 | 104 | 0.21 (0.07) |
| 3 | 72 | 42 | 0.12 (0.10) |
| 4 | 97 | 22 | 0.31 (0.11) |
| Average causal effect | | | 0.23 (0.05) |

Although in this particular application, the effect of CRM using a probit model is very close to the 23 percentage point increase estimated using the propensity score approach, there are two items worth noting. First, a probit analysis of $Y$ on $Z$ and $X$ provides a treatment effect that is conditional on the specific values of other covariates, while the propensity score method shown here provides an estimate of the average (over $X$) causal effect. If one were to compare the effect of CRM on $Y$ from both analyses on the probit scale and a linear additive model (in $X$) is appropriate, this distinction makes little difference. However, if the effect of CRM on $Y$ varied depending on whether a firm is in the service or manufacturing sector, for example, then a probit model with linear additive terms in $X$ is not appropriate. Even if this deficiency were corrected (say by including a CRM × MFG interaction term), the two analyses are not directly comparable, because they answer two different causal questions. Second, in general, the results of regression based models may not sustain "causal" interpretations, because they do not always ensure covariate balance across treatment and control firms (Dehejia and Wahba, 1999).

## 4. DISCUSSION

Our goal in this paper was to study the causal effect of CRM systems on one-to-one marketing effectiveness. We find evidence that CRM systems have indeed caused an improvement in one-to-one marketing effectiveness of firms. The use of potential outcomes reasoning and a propensity score approach offers four primary advantages to estimate the causal effect of sharp treatments such as CRM systems in a way that is not fully or satisfactorily covered by traditional approaches such as linear regression or probit models. First, the potential outcomes framework permits a more precise articulation of causal questions in terms of a comparison between two alternative states of the same firm (firm A with CRM versus firm A without CRM). As Rubin and Waterman (2006) elaborate, the causal effect is not a change in time of the observed outcome: instead, it is the difference between two potential outcomes, only one of which is observed. This clarifies why "before and after studies" (Connolly, 2003) and quasiexperimental studies [e.g., event studies that generate useful insights by studying changes in stock returns following an announcement regarding an electronic commerce application (Dehning, Richardson, Urbaczewski and Wells, 2004; Dehning, Richardson and Zmud, 2003; Im, Dow and Grover, 2001)] do not estimate "causal" effects in the sense described in this paper.

Second, although in the firm level application used in this paper we did not have many covariates, in general the propensity score reduces the dimensionality of observed covariates and avoids the problem of matching on multiple covariates. Because firms differ on multiple attributes, the propensity score approach offers a solution to the curse of dimensionality that has plagued much of the firm level research in the business value of IT literature, where researchers had to match on relatively few covariates to avoid losing degrees of freedom or statistical power.

Third, because we are primarily interested in estimating the causal effect of CRM systems, the propensity score approach allows us to focus on estimating that effect without having to specify how other observed covariates, such as IT investments and industry or other customer-facing IT systems that are correlated with CRM, may be related to the outcome variable one-to-one marketing effectiveness. In other words, the propensity score method allows researchers to escape strong functional form assumptions (e.g., constant treatment effect or linear effect of covariates) that are implicit in regression analysis (Dehejia and Wahba, 1999; Rubin, 1997).

Fourth and finally, the propensity score approach provides visibility of the extent to which CRM and non-CRM groups are similar to or different from each other based on observed covariates. Recall that while we started with a sample size of 487 firms, we had to drop six non-CRM firms because, for these firms, we did not find a matching CRM firm on support. We also ensured that CRM and non-CRM firms had a similar distribution of covariates within each stratum. In contrast, the extent to which CRM and non-CRM groups overlap is rarely, if ever, examined when analysts use classical regression methods.

Also note that the propensity score approach makes a separation between the balancing (or propensity score) stage, in which the goal is to ensure that treatment groups are comparable, and the analysis stage, in which the desired causal question is addressed. In contrast, with traditional regression methods, the researcher attempts to adjust simultaneously for selection bias and model the causal phenomena of interest. Unfortunately, the consideration of repeated models in this fashion (with the objective of controlling for covariates) may lead to unreliable treatment effect estimates. With the propensity score approach, on the other hand, a researcher is encouraged to "mine" the data (propose different propensity score models, consider different subclassifications of the propensity score, etc.) with the objective of achieving balance on observed pretreatment covariates. Since this first stage of analysis does not involve outcome data, there is less concern about inappropriately slanting the results of the analysis. Together, these advantages of the propensity score approach provide a powerful motivation for greater use of this approach to study causal effects not only in CRM and electronic commerce research, but also in other fields of management research.

Although this study extends previous research on the effect of CRM systems on firm performance (Boulding et al., 2005; Mithas, Krishnan and Fornell, 2005; Srinivasan and Moorman, 2005), several opportunities for future research remain. First, business value researchers have hinted at the possibility that, depending on their degree of preparedness, firms may benefit differentially from their IT systems (Lucas, 1993). Future research should systematically explore the issues related to treatment effect heterogeneity in the CRM context. Likewise, findings of observational studies, no matter how carefully done, are subject to the bias due to omission of unobservable variables and thus point to the need for a sensitivity analysis to assess the degree of potential bias. Our related work takes initial steps to address these issues of treatment effect heterogeneity and sensitivity analysis to assess the impact due to selection on unobservables at the firm and individual levels (Mithas and Krishnan, 2004a, b). Second, much of the CRM related work pertains to the business-to-consumer (B2C) domain and there are very few studies that have studied the phenomenon of CRM implementation in the business-to-business (B2B) domain (Mithas, Jones, Krishnan and Fornell, 2005). Since B2B transactions constitute a much greater proportion of economic activity compared with B2C transactions, it will be useful to understand the antecedents and consequences of CRM in the B2B context.

Going beyond the CRM context, we suggest that the potential outcome framework applies more broadly, from firm level questions to individual level and economy level causal questions such as those mentioned in the Introduction section. It is encouraging to note increasing use of the language of potential outcomes in academic research. For example, Bhagwati, Panagariya and Srinivasan (2004) seem to suggest that an economy level causal effect of offshoring on wages and jobs requires the use of the potential outcomes framework to pose the correct question. They note, "Forrester does not explain whether the prediction is that the U.S. economy will have 3.3 million fewer jobs in 2015 than it would otherwise have had because of outsourcing... or whether the prediction is that outsourcing will cause 3.3 million U.S. workers to shift from jobs that they might otherwise have had into different jobs..." (Bhagwati, Panagariya and Srinivasan, 2004, page 97). Similarly, at the individual level, researchers are interested in investigating the causal effect of the MBA degree for IT professionals. Use of the propensity score approach in this context can be informative because it shifts attention from estimation of causal effects based on a "before and after" type analysis to the potential outcome analysis used in this paper (Connolly, 2003; Mithas and Krishnan, 2004b; Pfeffer and Fong, 2003).

We acknowledge some limitations and challenges in the use of a propensity score approach in our CRM application. First, the most important limitation of this paper is that we have a relatively small set of observed covariates as is typical in firm level studies (Mithas, Krishnan and Fornell, 2005). Therefore, we are unable to take complete advantage of the power of the propensity score approach as a dimension reduction tool. Second, we focus on binary treatment and outcome variables in this application. However, the propensity score approach has been extended for ordinal, categorical and arbitrary treatments, and is also applicable to continuous outcomes. Third, as in other studies (Rosenbaum and Rubin, 1983a), we did not specifically account for the uncertainty in estimation of propensity scores in our estimation of standard errors of the overall causal effect. However, recent work suggests several methods to compute standard errors considering various sources of uncertainty, and finds that standard errors computed using different methods often give similar results (Agodini and Dynarski, 2004; Benjamin, 2003).

To conclude, this paper assessed the causal effect of CRM systems on one-to-one marketing effectiveness using a propensity score approach. We provided a detailed procedure to carry out the use of propensity score matching to assess the causal effect of CRM systems using a real data set from electronic commerce research. This paper illustrates the usefulness of the potential outcome approach and propensity score stratification to pose causal questions and to estimate causal effects of managerial and electronic commerce related IT interventions from observational data. We hope that this paper will encourage electronic commerce and management researchers to utilize this approach to answer their substantively interesting research questions in causal terms.

## ACKNOWLEDGMENTS

We thank Wolfgang Jank, Galit Shmueli and four anonymous reviewers for their guidance and help to improve this manuscript. We thank Jonathan Whitaker for helpful comments. We thank *InformationWeek* for providing the data for this research. In particular, we thank Rusty Weston, Bob Evans, Brian Gillooly, Stephanie Stahl, Lisa Smith and Helen D'Antoni for their help with data and for sharing their insights.

## REFERENCES

Agodini, R. and Dynarski, M. (2004). Are experiments the only option? A look at dropout prevention programs. *Rev. Econom. Statist.* **86** 180–194.

Angrist, J. D. and Krueger, A. B. (1999). Empirical strategies in labor economics. In *Handbook of Labor Economics* **3A** (O. Ashenfelter and D. Card, eds.) 1277–1366. North-Holland, Amsterdam.

Banker, R. D., Kauffman, R. J. and Mahmood, M. A. (1993). Measuring the business value of IT: A future oriented perspective. In *Strategic Information Technology Management*: *Perspectives on Organizational Growth and Competitive Advantage* (R. D. Banker, R. J. Kauffman and M. A. Mahmood, eds.) 595–605. Idea Group Publishing, Harrisburg, PA.

Bapna, R., Goes, P. and Gupta, A. (2001). Insights and analysis of online auctions. *Comm. ACM* **44**(11) 42–50.

Bapna, R., Jank, W. and Shmueli, G. (2005). Consumer surplus in online auctions. Working paper, Dept. Operations and Information Management, Univ. Connecticut. Available at www.sba.uconn.edu/users/rbapna/research.htm.

Barua, A. and Mukhopadhyay, T. (2000). Information technology and business performance: Past, present, and future. In *Framing the Domains of IT Management*: *Projecting the Future... Through the Past* (R. W. Zmud, ed.) 65–84. Pinnaflex, Cincinnati, OH.

Becker, S. O. and Ichino, A. (2002). Estimation of average treatment effects based on propensity scores. *The Stata Journal* **2** 358–377.

Benjamin, D. J. (2003). Does 401(k) eligibility increase saving? Evidence from propensity score subclassification. *J. Public Economics* **87** 1259–1290.

Bhagwati, J., Panagariya, A. and Srinivasan, T. N. (2004). The muddles over outsourcing. *J. Economic Perspectives* **18**(4) 93–114.

Boulding, W., Staelin, R., Ehret, M. and Johnston, W. J. (2005). A customer relationship management roadmap: What is known, potential pitfalls, and where to go. *J. Marketing* **69**(4) 155–166.

Brynjolfsson, E. and Hitt, L. M. (1998). Beyond the productivity paradox. *Comm. ACM* **41**(8) 49–55.

Connolly, M. (2003). The end of the MBA as we know it? *Academy of Management Learning and Education* **2** 365–367.

Dehejia, R. H. and Wahba, S. (1999). Causal effects in nonexperimental studies: Re-evaluating the evaluation of training programs. *J. Amer. Statist. Assoc.* **94** 1053–1062.

Dehejia, R. H. and Wahba, S. (2002). Propensity score-matching methods for nonexperimental causal studies. *Rev. Econom. Statist.* **84** 151–161.

Dehning, B., Richardson, V. J., Urbaczewski, A. and Wells, J. D. (2004). Reexamining the value relevance of e-commerce initiatives. *J. Management Inform. Syst.* **21** 55–82.

Dehning, B., Richardson, V. J. and Zmud, R. W. (2003). The value relevance of announcements of transformational information technology investments. *MIS Quarterly* **27** 637–656.

Fornell, C., Mithas, S., Morgeson, F. and Krishnan, M. S. (2006). Customer satisfaction and stock prices: High returns, low risk. *J. Marketing* **70**(1) 3–14.

Gregor, S. (2006). The nature of theory in information systems. *MIS Quarterly*. To appear.

Harvard Management Update (2000). A crash course in customer relationship management. *Harvard Management Update* **5**(3) 12 pages.

Heckman, J. J. (2005). The scientific model of causality (with discussion). *Sociological Methodology* **35** 1–162.

Holland, P. (1986). Statistics and causal inference (with discussion). *J. Amer. Statist. Assoc.* **81** 945–970. MR0867618

Im, K. S., Dow, K. E. and Grover, V. (2001). Research report: A reexamination of IT investment and the market value of the firm—An event study methodology. *Information Systems Research* **12** 103–117.

Imbens, G. W. (2004). Nonparametric estimation of average treatment effects under exogeneity: A review. *Rev. Econom. Statist.* **86** 4–29.

Jank, W. and Shmueli, G. (2006). Functional data analysis in electronic commerce research. *Statist. Sci.* **21** 155–166.

CRM SYSTEMS AND MARKETING EFFECTIVENESS 11

Kauffman, R. J. and Weill, P. (1989). An evaluative framework for research on the performance effects of information technology investments. In *Proc. International Conference on Information Systems*, Association for Information Systems 377–388. Boston.

Kohli, A. K. and Jaworski, B. J. (1990). Market orientation: The construct, research propositions, and managerial implications. *J. Marketing* **54**(2) 1–18.

Koppius, O. R. and Van Heck, E. (2002). Information architecture and electronic market performance in multidimensional auctions. Working paper, Erasmus Univ.

Lucas, H. C. (1975). Performance and the use of an information system. *Management Sci.* **21** 908–919.

Lucas, H. C. (1993). The business value of information technology: A historical perspective and thoughts for future research. In *Strategic Information Technology Management*: *Perspectives on Organizational Growth and Competitive Advantage* (R. D. Banker, R. J. Kauffman and M. A. Mahmood, eds.) 359–374. Idea Group Publishing, Harrisburg, PA.

Lucas, H. C. and Nielsen, N. R. (1980). The impact of the mode of information presentation on learning and performance. *Management Sci.* **26** 982–993.

Mithas, S. and Jones, J. L. (2006). Do auction parameters affect buyer surplus in e-auctions for procurement? *Production and Operations Management*. To appear.

Mithas, S., Jones, J. L., Krishnan, M. S. and Fornell, C. (2005). A theoretical integration of technology adoption and business value literature: The case of CRM systems. Working paper, Ross School of Business, Univ. Michigan.

Mithas, S. and Krishnan, M. S. (2004a). Causal effect of CRM systems on cross selling effectiveness and sales-force productivity by bounding a matching estimator. In *Proc. Ninth Annual INFORMS Conference on Information Systems and Technology* (H. Bhargava, C. Forman, R. Kauffman and D. J. Wu, eds.) 1–32. Denver, CO.

Mithas, S. and Krishnan, M. S. (2004b). Returns to managerial and technical competencies of information technology professionals: An empirical analysis. Working paper, Ross School of Business, Univ. Michigan.

Mithas, S., Krishnan, M. S. and Fornell, C. (2005). Why do customer relationship management applications affect customer satisfaction? *J. Marketing* **69**(4) 201–209.

Mithas, S. and Whitaker, J. (2006). Effect of information intensity and physical presence need on the global disaggregation of services: Theory and empirical evidence. Working paper, Smith School of Business, Univ. Maryland.

Peppers, D., Rogers, M. and Dorf, B. (1999). Is your company ready for one-to-one marketing? *Harvard Business Review* **77**(1) 3–12.

Pfeffer, J. and Fong, C. (2003). Assessing business schools: A reply to Connolly. *Academy of Management Learning and Education* **2** 368–370.

Rai, A., Patnayakuni, R. and Patnayakuni, N. (1997). Technology investment and business performance. *Comm. ACM* **40**(7) 89–97.

Rosenbaum, P. R. (1999). Choice as an alternative to control in observational studies (with discussion). *Statist. Sci.* **14** 259–304.

Rosenbaum, P. (2002). *Observational Studies*, 2nd ed. Springer, New York. MR1899138

Rosenbaum, P. R. and Rubin, D. B. (1983a). Assessing sensitivity to an unobserved binary covariate in an observational study with binary outcome. *J. Roy. Statist. Soc. Ser. B* **45** 212–218.

Rosenbaum, P. R. and Rubin, D. B. (1983b). The central role of the propensity score in observational studies for causal effects. *Biometrika* **70** 41–55. MR0742974

Rosenbaum, P. R. and Rubin, D. B. (1984). Reducing bias in observational studies using subclassification on the propensity score. *J. Amer. Statist. Assoc.* **79** 516–524.

Rubin, D. B. (1974). Estimating causal effects of treatments in randomized and nonrandomized studies. *J. Educational Psychology* **66** 688–701.

Rubin, D. B. (1977). Assignment to treatment group on the basis of a covariate. *J. Educational Statistics* **2** 1–26.

Rubin, D. B. (1997). Estimating causal effects from large data sets using propensity scores. *Ann. Internal Medicine* **127** 757–763.

Rubin, D. B. (2001). Using propensity scores to help design observational studies: Application to the tobacco litigation. *Health Services and Outcomes Research Methodology* **2** 169–188.

Rubin, D. B. (2005). Causal inference using potential outcomes: Design, modeling, decisions. *J. Amer. Statist. Assoc.* **100** 322–331. MR2166071

Rubin, D. B. and Waterman, R. P. (2006). Estimating the causal effects of marketing interventions using propensity score methodology. *Statist. Sci.* **21** 206–222.

Rust, R. T. and Kannan, P. K. (2003). E-service: A new paradigm for business in the electronic environment. *Comm. ACM* **46**(6) 36–42.

Sambamurthy, V., Bharadwaj, A. and Grover, V. (2003). Shaping agility through digital options: Reconceptualizing the role of information technology in contemporary firms. *MIS Quarterly* **27** 237–263.

Srinivasan, R. and Moorman, C. (2005). Strategic firm commitments and rewards for customer relationship management in online retailing. *J. Marketing* **69**(4) 193–200.

Venkatraman, N. (2004). Offshoring without guilt. *Sloan Management Review* **45**(3) 14–16.

Winship, C. and Morgan, S. L. (1999). The estimation of causal effects from observational data. *Annual Review of Sociology* **25** 659–706.